# Conformable Fractional Bessel Equation and Bessel Functions

Ahmet Gökdoğan[a], Emrah Ünal[b], Ercan Çelik[c]

[a] Department of Mathematical Engineering, Gümüşhane University, 29100 Gümüşhane, Turkey,
gokdogan@gumushane.edu.tr
[b] Department of Elementary Mathematics Education, Artvin Çoruh University, 08100 Artvin, Turkey
emrah.unal@artvin.edu.tr
[c] Department of Mathematics, Atatürk University, 25400 Erzurum, Turkey,
ecelik@atauni.edu.tr

**Abstract**

In this work, we study the fractional power series solutions around regular singular point $x = 0$ of conformable fractional Bessel differential equation and fractional Bessel functions. Then, we compare fractional solutions with ordinary solutions. In addition, we present certain property of fractional Bessel functions.

**Keywords:** conformable fractional derivative, conformable fractional Bessel equation, conformable fractional Bessel functions

## 1. Introduction

The idea of fractional derivatives was raised first by L'Hospital in 1695. The fractional calculus and its mathematical consequences attracted many mathematicians such as Fourier, Euler, Laplace. Various definitions of non-integer order integral or derivative was given by many mathematicians [1]. For fractional derivative, the most popular definitions are Riemann-Liouville and Caputo definitions as follows, respectively:

(I) Riemann Liouville definition:
$$D_x^\alpha f(x) = \frac{1}{\Gamma(n-\alpha)} \left(\frac{d}{dx}\right)^n \int_0^x (x-t)^{n-\alpha-1} f(t) dt, \quad n-1 < \alpha \leq n;$$

(II) Caputo definition:
$$D_x^\alpha f(x) = \frac{1}{\Gamma(n-\alpha)} \int_0^x (x-t)^{n-\alpha-1} \left(\frac{d}{dx}\right)^n f(t) dt, \quad n-1 < \alpha \leq n.$$

For Riemann-Liouville, Caputo and other definitions and the characteristics of these definitions, we refer to reader to [2-4]

Recently, Khalil and colleagues have given a new definition of non-integer order integral and derivative [5]. This new definition is used a limit form as in usual derivatives. They also proved the product rule, the fractional Rolle theorem and mean value theorem. This new theory is improvised by Abdeljawad. For instance, he provides definitions of left and right conformable fractional derivatives and fractional integrals of higher order (i.e. of order $\alpha > 1$), Taylor power series representation and Laplace transform of few certain functions, fractional integration by parts formulas, chain rule and Gronwall inequality in [6].

In [7], authors have given existence and uniqueness theorems of sequential linear conformable fractional differential equations. The existence of fractional power series solutions around an ordinary point of sequential conformable fractional differential equation of order $2\alpha$, the solution of sequential conformable fractional Hermite differential equation and fractional

Hermite polynomials are analyzed in [8]. Legendre conformable fractional equation and Legendre fractional polynomials are studied in [9].

In [10], the existence of fractional power series solutions around a regular-singular point of sequential conformable fractional differential equation of order $2\alpha$ is analyzed. In this paper, we give solutions of the conformable fractional Bessel equation and certain properties of conformable fractional Bessel functions.

## 2. Conformable Fractional Calculus

**Definition 2.1.** [5] Let a function $f: [a, \infty) \to \mathbb{R}$ be given. Then, the left conformable fractional derivative of $f$ order $\alpha$ is defined by

$$(T_\alpha^a f)(x) = \lim_{\varepsilon \to 0} \frac{f(x + \varepsilon(x-a)^{1-\alpha}) - f(x)}{\varepsilon}$$

for all $x > a, \alpha \in (0,1]$. When $a = 0$, it is written as $T_\alpha$. If $(T_\alpha f)(x)$ exists on $(a, b)$, then $(T_\alpha^a f)(a) = \lim_{x \to a^+} (T_\alpha^a f)(x)$.

**Definition 2.2.** [5] Given a function $f: (-\infty, b] \to \mathbb{R}$. Then the right conformable fractional derivative of $f$ order $\alpha$ is defined by

$$(_\alpha^b T f)(x) = -\lim_{\varepsilon \to 0} \frac{f(x + \varepsilon(b-x)^{1-\alpha}) - f(x)}{\varepsilon}$$

for all $x < b, \alpha \in (0,1]$. If $(_\alpha T f)(x)$ exists on $(a, b)$, then $(_\alpha^b T f)(b) = \lim_{x \to b^-} (_\alpha^b T f)(x)$.

**Theorem 2.1.** [4] Let $\alpha \in (0,1]$ and $f, g$ be $\alpha$-differentiable at a point $x > 0$. Then

(1) $\frac{d^\alpha}{dx^\alpha}(af + bg) = a\frac{d^\alpha f}{dx^\alpha} + b\frac{d^\alpha g}{dx^\alpha}$, for all $a, b \in \mathbb{R}$

(2) $\frac{d^\alpha}{dx^\alpha}(x^p) = px^{p-\alpha}$, for all $p \in \mathbb{R}$

(3) $\frac{d^\alpha}{dx^\alpha}(\lambda) = 0$, for all constant functions $f(x) = \lambda$

(4) $\frac{d^\alpha}{dx^\alpha}(fg) = f\frac{d^\alpha}{dx^\alpha}(g) + g\frac{d^\alpha}{dx^\alpha}(f)$

(5) $\frac{d^\alpha}{dx^\alpha}(f/g) = \frac{g\frac{d^\alpha}{dx^\alpha}(f) - f\frac{d^\alpha}{dx^\alpha}(g)}{g^2}$

(6) If, in addition, $f$ is differentiable, then $\frac{d^\alpha}{dx^\alpha}(f(x)) = x^{1-\alpha}\frac{df}{dx}(x)$.

**Theorem 2.2.** [5] Assume that $f$ is infinitely $\alpha$-differentiable function, for some $0 < \alpha \leq 1$ at a neighborhood of a point $x_0$. Then, $f$ has the fractional power series expansion:

$$f(x) = \sum_{k=0}^{\infty} \frac{\left(^{(k)}T_\alpha^{x_0}f\right)(x_0)(x-x_0)^{k\alpha}}{\alpha^k k!}, \qquad x_0 < x < x_0 + R^{1/\alpha}, \qquad R > 0.$$

Here, $\left(^{(k)}T_\alpha^{x_0}f\right)(x_0)$ means the application of the fractional derivative $k$ times.

## 3. Conformable Fractional Bessel Equation and Conformable Fractional Bessel Functions

Consider the sequential conformable fractional Bessel equation

$$x^{2\alpha}T_\alpha T_\alpha y + \alpha x^\alpha T_\alpha y + \alpha^2(x^{2\alpha} - p^2)y = 0, \tag{1}$$

where $0 < \alpha \leq 1$ and $p$ is any real number. If $\alpha = 1$, then the equation (1) is classical Bessel equation.

$x = 0$ is a $\alpha$ – regular singular point for the equation. In this case, for $x > 0$, we can investigate solution by a fractional Frobenius series as follows:

$$y = x^{r\alpha} \sum_{n=0}^{\infty} c_n x^{n\alpha} = \sum_{n=0}^{\infty} c_n x^{(n+r)\alpha}.$$

We let

$$T_\alpha y = \sum_{n=0}^{\infty} \alpha(n+r)c_n x^{(n+r-1)\alpha},$$

$$T_\alpha T_\alpha y = \sum_{n=0}^{\infty} \alpha^2(n+r)(n+r-1)c_n x^{(n+r-2)\alpha}.$$

Now, substituting these expressions in equation (1), we get

$$\sum_{n=0}^{\infty} \alpha^2(n+r)(n+r-1)c_n x^{(n+r)\alpha} + \sum_{n=0}^{\infty} \alpha^2(n+r)c_n x^{(n+r)\alpha} + \sum_{n=0}^{\infty} \alpha^2 c_n x^{(n+r+2)\alpha}$$

$$- \sum_{n=0}^{\infty} \alpha^2 p^2 c_n x^{(n+r)\alpha} = 0$$

$$\Rightarrow \sum_{n=0}^{\infty} \alpha^2(n+r)(n+r-1)c_n x^{(n+r)\alpha} + \sum_{n=0}^{\infty} \alpha^2(n+r)c_n x^{(n+r)\alpha} + \sum_{n=2}^{\infty} \alpha^2 c_{n-2} x^{(n+r)\alpha}$$

$$- \sum_{n=0}^{\infty} \alpha^2 p^2 c_n x^{(n+r)\alpha} = 0$$

$$\Rightarrow (r(r-1)\alpha^2 + r\alpha^2 - \alpha^2 p^2)c_0 x^{r\alpha} + (r(r+1)\alpha^2 + (r+1)\alpha^2 - \alpha^2 p^2)c_1 x^{(r+1)\alpha}$$

$$+ \sum_{n=0}^{\infty} [(\alpha^2(n+r)(n+r-1) + (n+r)\alpha^2 - \alpha^2 p^2)c_n + \alpha^2 c_{n-2}]x^{(n+r)\alpha}$$

$$= 0.$$

If we get
$$I(r) = r(r-1)\alpha^2 + r\alpha^2 - \alpha^2 p^2, \qquad (2)$$
Then, the last equation can be rewritten as follow:
$$I(r)c_0 x^{r\alpha} + I(r+1)c_1 x^{(r+1)\alpha} + \sum_{n=0}^{\infty}[I(r+n)c_n + \alpha^2 c_{n-2}]x^{(n+r)\alpha} = 0.$$

Let us $c_0 \neq 0$, then we obtain
$$I(r) = 0.$$
Since it is that $\alpha^2 \neq 0$, we can write
$$I(r) = 0 \Rightarrow r(r-1)\alpha^2 + r\alpha^2 - \alpha^2 p^2 = 0 \Rightarrow r^2 - p^2 = 0.$$
Hence, we find
$$r_1 = p, r_2 = -p.$$
Firstly, let us analyze $p = 0$ case and find solutions of fractional Bessel Equation of order zero. In this case, roots of fractional indicial equation are
$$r_1 = 0, r_2 = 0.$$
For $r_1 = 0$, we obtain
$$I(r_1 + 1)c_1 = 0 \Rightarrow c_1 = 0.$$
In addition to this, the following recurrence relation is valid
$$c_n = -\frac{\alpha^2 c_{n-2}}{I(r_1+n)} = -\frac{c_{n-2}}{(n)^2}, \quad n \geq 2.$$
By using the recurrence relation, we have
$$c_3 = c_5 = \cdots = 0.$$
Since $c_0$ is arbitrary constant, we get
$$c_2 = -\frac{c_0}{2^2},$$
$$c_4 = -\frac{c_2}{4^2} = \frac{c_0}{2^4 2! \, 2!},$$
$$c_6 = -\frac{c_4}{6^2} = -\frac{c_0}{2^6 3! \, 3!},$$
$$\vdots$$
$$c_{2n} = \frac{c_0(-1)^n}{2^{2n} n! \, n!}.$$
Hence, the first solution of fractional Bessel Equation of order zero is
$$y_1 = c_0 \sum_{n=0}^{\infty} \frac{(-1)^n x^{2n\alpha}}{2^{2n} n! \, n!} = c_0 \sum_{n=0}^{\infty} \frac{(-1)^n}{n! \, n!} \left(\frac{x^\alpha}{2}\right)^{2n},$$
and the Bessel functions of order zero is found by
$$(J_\alpha)_0(x) = \sum_{n=0}^{\infty} \frac{(-1)^n}{n! n!} \left(\frac{x^\alpha}{2}\right)^{2n}.$$

For the second solution, we use the general form that [10]. The general form has the next form
$$y_2(x; s_1) = \ln|x| y_1(x; s_1) + \sum_{n=0}^{\infty} b_n(x)^{(n+s_1+1)\alpha}. \tag{3}$$
For $s_1 = 0$, applying conformable fractional derivatives to equation (3), we have
$$T_\alpha y_2 = T_\alpha y_1 \ln x + \frac{y_1}{x^\alpha} + \sum_{n=0}^{\infty} \alpha(n+1) b_n x^{n\alpha},$$
$$T_\alpha T_\alpha y_2 = T_\alpha T_\alpha y_1 \ln x + 2\frac{T_\alpha y_1}{x^\alpha} - \alpha \frac{y_1}{x^{2\alpha}} + \sum_{n=1}^{\infty} \alpha^2 n(n+1) b_n x^{(n-1)\alpha}.$$
Substituting these expressions in equation (1), we get
$$2x^\alpha \sum_{n=1}^{\infty} \frac{2n\alpha(-1)^n}{n!n!2^{2n}} x^{(2n-1)\alpha} + \sum_{n=0}^{\infty} \alpha^2(n+1) b_n x^{(n+1)\alpha} + \sum_{n=1}^{\infty} \alpha^2 n(n+1) b_n x^{(n+1)\alpha} + \sum_{n=0}^{\infty} \alpha^2 b_n x^{(n+3)\alpha} = 0.$$
If we make arrangements to the expression above, then we get
$$\alpha^2 b_0 x^\alpha + 2^2 \alpha^2 b_1 x^{2\alpha} - \alpha x^{2\alpha} + \sum_{n=2}^{\infty} \frac{4n\alpha(-1)^n}{n!n!2^{2n}} x^{2n\alpha} + \sum_{n=2}^{\infty} (\alpha^2(n+1)^2 b_n + \alpha^2 b_{n-2}) x^{(n+1)\alpha} = 0.$$
From this equation, $b_0 = 0$ and $b_1 = \frac{1}{\alpha 2^2}$ are obtained. Besides, if $n$ is even, then $\alpha^2(n+1)^2 b_n + \alpha^2 b_{n-2} = 0$. This implies that $b_0 = b_2 = \cdots = 0$. If $n$ is odd, for $n \geq 2$, then
$$\alpha^2 (2n)^2 b_{2n-1} + \alpha^2 b_{2n-3} = \frac{(-1)^{n+1} 4n\alpha}{n!n!2^{2n}}.$$
Hence, we get
$$b_3 = -\frac{1}{\alpha 2^2 4^2}\left(1 + \frac{1}{2}\right),$$
$$b_5 = \frac{1}{\alpha 2^2 4^2 6^2}\left(1 + \frac{1}{2} + \frac{1}{3}\right),$$
$$\vdots$$
$$b_{2n-1} = \frac{(-1)^{n+1}}{\alpha 2^{2n} n! n!}\left(1 + \frac{1}{2} + \cdots + \frac{1}{n}\right).$$
Let us $\left(1 + \frac{1}{2} + \cdots + \frac{1}{n}\right) = H_n$. Then, the second solution of Bessel Equation of order zero is
$$y_2(x) = (J_\alpha)_0(x) \ln x + \frac{1}{\alpha} \sum_{n=1}^{\infty} \frac{(-1)^{n+1} H_n}{2^{2n} n! n!} x^{2n\alpha}.$$
**Corollary 3.1.** Let $Y_2(x)$ be the second solution for classical Bessel equation of order zero, then
$$y_2(x) = \frac{1}{\alpha} Y_2(x^\alpha).$$
Now, for $p \neq 0$, let us analyze solutions of fractional Bessel Equation of order $p$. In this case for $p > 0$,
$$r_1 = p, r_2 = -p, \quad r_1, r_2 \in \mathbb{R}.$$
For $r_1 = p$, we have
$$I(r_1 + 1)c_1 = [\alpha^2(p)(p+1) + (p+1)\alpha^2 - \alpha^2 p^2]c_1 = 0$$
$$(2p+1)c_1 = 0.$$
Because of $p > 0$, it follows that $c_1 = 0$. The recurrence relation is
$$c_n = -\frac{\alpha^2 c_{n-2}}{I(r+n)} = -\frac{c_{n-2}}{n(n+2p)}.$$

From $c_1 = 0$ and last recurrence relation, we get
$$c_3 = c_5 = \cdots = 0,$$
and
$$c_2 = -\frac{c_0}{2(2+2p)} = -\frac{c_0}{2^2 1!(p+1)},$$
$$c_4 = -\frac{c_2}{4(4+2p)} = \frac{c_0}{2^4 2!(p+1)(p+2)},$$
$$c_6 = -\frac{c_4}{6(6+2p)} = -\frac{c_0}{2^6 3!(p+1)(p+2)(p+3)},$$
$$\vdots$$
$$c_{2n} = \frac{(-1)^n c_0}{2^{2n} n!(p+1)(p+2)\ldots(p+n)}.$$

Let us $c_0 = \frac{c}{2^p \Gamma(p+1)}$. Hence, the first solution of fractional Bessel equation of order $p$ has the form

$$y_1 = cx^{p\alpha} \sum_{n=0}^{\infty} \frac{(-1)^n x^{2n\alpha}}{2^{2n+p} n! \Gamma(p+1)(p+1)(p+2)\ldots(p+n)} = c \sum_{n=0}^{\infty} \frac{(-1)^n}{n! \Gamma(p+n+1)} \left(\frac{x^\alpha}{2}\right)^{2n+p}.$$

Besides, Bessel functions of order $p$ is valid

$$(J_\alpha)_p(x) = \sum_{n=0}^{\infty} \frac{(-1)^n}{n! \Gamma(p+n+1)} \left(\frac{x^\alpha}{2}\right)^{2n+p}. \tag{4}$$

For $r_2 = -p$, if $r_1 - r_2 = 2p$ is not a positive integer, then we have
$$I(r_2+1)c_1 = [\alpha^2(-p)(-p+1) + (-p+1)\alpha^2 - \alpha^2 p^2]c_1 = 0$$
$$(1-2p)c_1 = 0.$$

This implies that $c_1 = 0$. For $n \geq 2$, the recurrence relation is
$$c_n = -\frac{\alpha^2 c_{n-2}}{I(r+n)} = -\frac{c_{n-2}}{n(n-2p)}.$$

Because of $c_1 = 0$, we can write
$$c_3 = c_5 = \cdots = 0.$$

If $n$ is even, then the coefficient is found as
$$c_2 = -\frac{2c_0}{2^2 1!(-p+1)},$$
$$c_4 = -\frac{c_2}{4(4-2p)} = \frac{c_0}{2^4 2!(-p+1)(-p+2)},$$
$$c_6 = -\frac{c_4}{6(6-2p)} = -\frac{c_0}{2^6 3!(-p+1)(-p+2)(-p+3)},$$
$$\vdots$$
$$c_{2n} = \frac{(-1)^n c_0}{2^{2n} n!(-p+1)(-p+2)(-p+3)\ldots(-p+n)}.$$

Let us $c_0 = \frac{c}{2^{-p} \Gamma(-p+1)}$. In this case, the second solution of fractional Bessel equation of order $p$ is

$$y_2 = cx^{-p\alpha} \sum_{n=0}^{\infty} \frac{(-1)^n x^{2n\alpha}}{2^{2n-p} n! \Gamma(-p+1)(-p+1)(-p+2)(-p+3)\ldots(-p+n)} = c \sum_{n=0}^{\infty} \frac{(-1)^n}{n! \Gamma(-p+n+1)} \left(\frac{x^\alpha}{2}\right)^{2n+p}.$$

The second type Bessel functions of order $p$ has the form

$$(J_\alpha)_{-p}(x) = \sum_{n=0}^{\infty} \frac{(-1)^n}{n!\,\Gamma(-p+n+1)} \left(\frac{x^\alpha}{2}\right)^{2n-p}. \tag{5}$$

Since $\Gamma(-p + n + 1)$ exists for $n \geq 0$, provided $p$ is not a positive integer, we see that $(J_\alpha)_{-p}$ exists in this case, even though $r_1 - r_2 = 2p$ is a positive integer. At that, the second solution also exists as follows

$$(J_\alpha)_{-p}(x) = \sum_{n=0}^{\infty} \frac{(-1)^n}{n!\,\Gamma(-p+n+1)} \left(\frac{x^\alpha}{2}\right)^{2n-p}.$$

Hence the last remaining case is that for which $p$ is a positive integer. In this case, we use the form in [10]. That is,

$$y_2(x; r_2) = C.\ln|x| y_1(x; r_1) + \sum_{n=0}^{\infty} b_n(x)^{(n-m)\alpha} \tag{6}$$

Let us $p = m$. Hence

$$(J_\alpha)_m(x) = \sum_{n=0}^{\infty} c_{2n}(x^\alpha)^{2n+m},$$

where

$$c_{2n} = \frac{(-1)^n}{2^{2n+m} n!\,(m+n)!}. \tag{7}$$

Substituting the equation (6) and its conformable fractional derivatives to the equation (1), we get

$$\sum_{n=0}^{\infty} \frac{2C(2n+m)\alpha(-1)^n}{2^{2n+m} n!\,(m+n)!} x^{2(n+m)\alpha} + \sum_{n=2}^{\infty}[\alpha^2 n(n-2m)b_n + \alpha^2 b_{n-2}]x^{n\alpha} + \alpha^2(1-2m)b_1 x^\alpha = 0. \tag{8}$$

$b_1 = 0$ is found by the last equation. In addition, if $m > 1$, we have

$$\alpha^2 n(n - 2m)b_n + \alpha^2 b_{n-2} = 0 \quad (n = 2,3,\ldots,2m-1).$$

This implies that

$$b_1 = b_3 = \cdots = b_{2m-1} = 0.$$

Besides, If $n$ is even, then the remaining coefficient is finded as:

$$b_2 = \frac{b_0}{2^2(m-1)},$$

$$b_4 = \frac{b_0}{2^4 2!\,(m-1)(m-2)},$$

$$\vdots$$

$$b_{2j} = \frac{b_0}{2^{2j} j!\,(m-1)(m-2)\ldots(m-j)}.$$

where $j = 1,2,\ldots,m-1$. Now, for the equation (8), let us compare coefficients of $x^{2m\alpha}$. Substituting $n = 0$ in the first series of equation (8) and $n = 2m$ in the second series of equation (8), we obtain

$$(\alpha^2 b_{2m-2} + 2C\alpha m c_0)x^{2m\alpha} = 0.$$

That is,

$$\alpha b_{2m-2} = -\frac{C}{2^{m-1}(m-1)!}.$$

For $j = m - 1$,
$$\alpha \frac{b_0}{2^{2m-2}(m-1)!(m-1)!} = -\frac{C}{2^{m-1}(m-1)!}.$$

Hence, we have
$$C = -\frac{\alpha b_0}{2^{m-1}(m-1)!}. \qquad (9)$$

Since first series in equation (8) contains only even power of $x^\alpha$, we have
$$b_{2m+1} = b_{2m+3} = \cdots = 0.$$

The coefficient $b_{2m}$ is undetermined, but the remaining coefficients $b_{2m+2}, b_{2m+4}, \ldots$ are defined by following equation:
$$\alpha^2 2n(2n + 2m)b_{2n+2m} + \alpha^2 b_{2n+2m-2} = -2C(2n + m)\alpha c_{2n}, \qquad (10)$$

where $= 1, 2, \ldots$.

For $n = 1$, from equation (10), we obtain
$$b_{2m+2} = -\frac{C}{2\alpha} c_2 \left(1 + \frac{1}{m+1}\right) - \frac{b_{2m}}{4(1+m)}.$$

Now, let us choose $b_{2m}$ so that
$$\frac{b_{2m}}{4(1+m)} = \frac{C}{2\alpha} c_2 \left(1 + \frac{1}{2} + \cdots + \frac{1}{m}\right). \qquad (11)$$

From equation (7), $4(1 + m)c_2 = -c_0$ is obtained. Using this expression in equation (11), we get
$$b_{2m} = -\frac{C}{2\alpha} c_0 \left(1 + \frac{1}{2} + \cdots + \frac{1}{m}\right).$$

With this choose of $b_{2m}$, we have
$$b_{2m+2} = -\frac{C}{2\alpha} c_2 \left(1 + 1 + \frac{1}{2} + \cdots + \frac{1}{m+1}\right).$$

For $n = 2$, we obtain
$$b_{2m+4} = -\frac{C}{2\alpha} c_4 \left(\frac{1}{2} + \frac{1}{m+2}\right) - \frac{b_{2m+2}}{2^2 2(2+m)}. \qquad (12)$$

From equation (7), $2^2 2(2 + m)c_4 = -c_2$ is obtained. Hence, we get
$$\frac{b_{2m+2}}{2^2 2(2+m)} = \frac{C}{2\alpha} c_4 \left(1 + 1 + \frac{1}{2} + \cdots + \frac{1}{m+1}\right). \qquad (13)$$

Substituting the equation (13) in the equation (12), we have
$$b_{2m+4} = -\frac{C}{2\alpha} c_4 \left(1 + \frac{1}{2} + 1 + \frac{1}{2} + \cdots + \frac{1}{m+2}\right).$$

And so on. It can be shown by induction that
$$b_{2m+2n} = -\frac{C}{2\alpha} c_{2n} \left(1 + \frac{1}{2} + \cdots + \frac{1}{m}\right)\left(1 + \frac{1}{2} + \cdots + \frac{1}{m+n}\right).$$

Hence, the second solution of the equation (1) is defined by

$$y_2(x) = b_0 x^{-m\alpha} + b_0 x^{-m\alpha} \sum_{j=1}^{m-1} \frac{x^{2j\alpha}}{2^{2j} j!(m-1)(m-2)\ldots(m-j)} - \frac{C}{2\alpha} c_0 \left(1 + \frac{1}{2} + \cdots + \frac{1}{m}\right) x^{m\alpha} - \frac{C}{2\alpha} \sum_{n=1}^{\infty} c_{2n} \left(1 + \frac{1}{2} + \cdots + \frac{1}{m}\right)\left(1 + \frac{1}{2} + \cdots + \frac{1}{m+n}\right) x^{(2n+m)\alpha} + c \ln x (J_\alpha)_m.$$

where $b_0$ and $C$ are constants which related by (9) and $c_{2n}$ is given by (7). If $C = 1$, then the second solution is indicated by $(K_\alpha)_m$ and, hence, we obtain

$$(K_\alpha)_m = -\frac{1}{2\alpha}\left(\frac{x^\alpha}{2}\right)^{-m} \sum_{j=0}^{m-1} \frac{(m-j-1)!}{j!} \left(\frac{x^\alpha}{2}\right)^{2j} - \frac{1}{2\alpha} \frac{1}{m!} c_0 \left(1 + \frac{1}{2} + \cdots + \frac{1}{m}\right)\left(\frac{x^\alpha}{2}\right)^m - \frac{1}{2\alpha}\left(\frac{x^\alpha}{2}\right)^m \sum_{n=1}^{\infty} \frac{(-1)^n}{n!(m+n)!}\left(1 + \frac{1}{2} + \cdots + \frac{1}{m}\right)\left(1 + \frac{1}{2} + \cdots + \frac{1}{m+n}\right)\left(\frac{x^\alpha}{2}\right)^{2n} + \ln x (J_\alpha)_m.$$

**Corollary 3.2.** Let $K_m(x)$ be the second solution for classical Bessel equation of order $m$, then, we can write

$$(K_\alpha)_m(x) = \frac{1}{\alpha} K_m(x^\alpha).$$

**Example 3.1.** Consider fractional Bessel functions of order $p = \pm\frac{1}{2}$. Show that

$$(J_\alpha)_{1/2}(x) = \sqrt{\frac{2}{\pi x^\alpha}} \sin(x^\alpha) \text{ and } (J_\alpha)_{-1/2}(x) = \sqrt{\frac{2}{\pi x^\alpha}} \cos(x^\alpha).$$

**Solution:** Substituting $p = \frac{1}{2}$ in (1), we get

$$(J_\alpha)_{1/2}(x) = \sum_{n=0}^{\infty} \frac{(-1)^n}{n! \Gamma\left(\frac{1}{2}+n+1\right)} \left(\frac{x^\alpha}{2}\right)^{2n+\frac{1}{2}}.$$

To simplify this expression, we use the following fact

$$\Gamma\left(\frac{1}{2} + n + 1\right) = \frac{(2n+1)!}{2^{2n+1} n!} \sqrt{\pi}.$$

Thus,

$$(J_\alpha)_{1/2}(x) = \frac{1}{\sqrt{\pi}} \sum_{n=0}^{\infty} \frac{(-1)^n 2^{2n+1}}{(2n+1)!} \left(\frac{x^\alpha}{2}\right)^{2n+\frac{1}{2}},$$

$$(J_\alpha)_{1/2}(x) = \sqrt{\frac{2}{\pi x^\alpha}} \sum_{n=0}^{\infty} \frac{(-1)^n}{(2n+1)!} (x^\alpha)^{(2n+1)},$$

$$(J_\alpha)_{1/2}(x) = \sqrt{\frac{2}{\pi x^\alpha}} \sin(x^\alpha).$$

The second relation is proved similarly by substituting $p = -\frac{1}{2}$ in (5).

**Property 3.1.** Let $p$ is a nonnegative integer. Then,

i. $T_\alpha\left(x^{p\alpha}(J_\alpha)_p(x)\right) = \alpha x^{p\alpha}(J_\alpha)_{p-1}(x)$

ii. $T_\alpha\left(x^{-p\alpha}(J_\alpha)_p(x)\right) = -\alpha x^{-p\alpha}(J_\alpha)_{p+1}(x)$

iii. $T_\alpha\left((J_\alpha)_p(x)\right) = \alpha(J_\alpha)_{p-1}(x) - \frac{\alpha p}{x^\alpha}(J_\alpha)_p(x)$

iv. $T_\alpha\left((J_\alpha)_p(x)\right) = \frac{\alpha p}{x^\alpha}(J_\alpha)_p(x) - \alpha(J_\alpha)_{p+1}(x)$

v. $(J_\alpha)_{p+1}(x) = \frac{2p}{x^\alpha}(J_\alpha)_p(x) - (J_\alpha)_{p-1}(x)$

vi. $(J_\alpha)_{-p}(x) = (-1)^p (J_\alpha)_p(x)$

**Proof:**

i. $x^{p\alpha}(J_\alpha)_p(x) = \sum_{n=0}^{\infty} \frac{(-1)^n (x^\alpha)^{2n+2p}}{n!(p+n)! 2^{2n+p}}$,

$$T_\alpha\left(x^{p\alpha}(J_\alpha)_p(x)\right) = \alpha \sum_{n=0}^{\infty} \frac{(-1)^n (x^\alpha)^{2n+2p-1}}{n!(p+n-1)! 2^{2n+p-1}}$$

$$= \alpha x^{p\alpha} \sum_{n=0}^{\infty} \frac{(-1)^n}{n! \Gamma(p+n)} \left(\frac{x^\alpha}{2}\right)^{2n+p-1}$$

$$= \alpha x^{p\alpha} (J_\alpha)_{p-1}(x).$$

ii. $x^{-p\alpha}(J_\alpha)_p(x) = \sum_{n=0}^{\infty} \frac{(-1)^n (x^\alpha)^{2n}}{n!(p+n)! 2^{2n+p}}$

$$T_\alpha\left(x^{-p\alpha}(J_\alpha)_p(x)\right) = \alpha \sum_{n=1}^{\infty} \frac{(-1)^n (x^\alpha)^{2n-1}}{(n-1)!(p+n)! 2^{2n+p-1}}$$

$$= -\alpha x^{-p\alpha} \sum_{n=0}^{\infty} \frac{(-1)^n}{n! \Gamma(p+n+2)} \left(\frac{x^\alpha}{2}\right)^{2n+p+1}$$

$$= -\alpha x^{-p\alpha} (J_\alpha)_{p+1}(x).$$

iii. Using Theorem 2.1. (4), we have
$$T_\alpha\left(x^{p\alpha}(J_\alpha)_p(x)\right) = p\alpha x^{(p-1)\alpha}(J_\alpha)_p(x) + x^{p\alpha} T_\alpha\left((J_\alpha)_p(x)\right).$$

From property 4.1. (i), we get
$$p\alpha x^{(p-1)\alpha}(J_\alpha)_p(x) + x^{p\alpha} T_\alpha\left((J_\alpha)_p(x)\right) = \alpha x^{p\alpha}(J_\alpha)_{p-1}(x).$$

Hence, we obtain
$$T_\alpha\left((J_\alpha)_p(x)\right) = \alpha (J_\alpha)_{p-1}(x) - \frac{\alpha p}{x^\alpha}(J_\alpha)_p(x).$$

iv. Similarly, using Theorem 2.1. (4), we have
$$T_\alpha\left(x^{-p\alpha}(J_\alpha)_p(x)\right) = -p\alpha x^{(-p-1)\alpha}(J_\alpha)_p(x) + x^{-p\alpha} T_\alpha\left((J_\alpha)_p(x)\right).$$

From property 4.1. (ii), we get
$$-p\alpha x^{(-p-1)\alpha}(J_\alpha)_p(x) + x^{-p\alpha} T_\alpha\left((J_\alpha)_p(x)\right) = -\alpha x^{-p\alpha}(J_\alpha)_{p+1}(x).$$

Hence, we obtain
$$T_\alpha\left((J_\alpha)_p(x)\right) = \frac{\alpha p}{x^\alpha}(J_\alpha)_p(x) - \alpha (J_\alpha)_{p+1}(x).$$

v. Subtracting from (iv) to (iii), we obtain the requested result.

vi. If $s$ is not positive integer, we have
$$(J_\alpha)_{-s}(x) = \sum_{n=0}^{\infty} \frac{(-1)^n}{n! \, \Gamma(-s+n+1)} \left(\frac{x^\alpha}{2}\right)^{2n-s}.$$

Recall that $\lim_{z \to -p} |\Gamma(z)| = \infty$ for $p = 0$ or $p \in Z^+$. When $s \to p$, $(-s + n + 1)$ tends to 0 or a negative integer for $n = 0,1,2, \ldots, (p-1)$. For such values of $n$, the coefficients of $\left(\frac{x^\alpha}{2}\right)^{2n-s}$ in the series above approaches 0:

$$\lim_{s \to p} \frac{(-1)^n}{n!\, \Gamma(-s+n+1)} = 0.$$

Then, we get

$$(J_\alpha)_{-p}(x) = \lim_{s \to p}(J_\alpha)_{-s}(x) = \sum_{n=p}^{\infty} \frac{(-1)^n}{n!\, \Gamma(-p+n+1)} \left(\frac{x^\alpha}{2}\right)^{2n-p}.$$

Now, replacing $n$ by $n + p$ in the above series and using the fundamental property of the Gamma function, we obtain

$$(J_\alpha)_{-p}(x) = (-1)^p \sum_{n=0}^{\infty} \frac{(-1)^n}{n!\, \Gamma(p+n+1)} \left(\frac{x^\alpha}{2}\right)^{2n+p} = (-1)^p (J_\alpha)_p(x).$$

## 6. Conclusion

In this work, we give power series solutions around regular-singular point $x = 0$ of the conformable fractional Bessel differential equation. In addition to this, we present certain properties of the conformable fractional Bessel functions. The findings of this study indicate that the results obtained in fractional case conform with results obtained in ordinary case.